\documentclass[oneside,12pt]{amsart}
\usepackage[a4paper]{geometry}     
\usepackage{amsfonts,amsmath,latexsym,amssymb}
\usepackage{amsthm,mathrsfs,upref,mathptmx}
\usepackage{enumerate,ifthen,cite}
\usepackage[mathscr]{eucal}
\usepackage[T1]{fontenc}
\usepackage{todonotes}

\newcommand{\N}{\mathbb{N}}
\newcommand{\Z}{\mathbb{Z}}
\newcommand{\Q}{\mathbb{Q}}
\newcommand{\R}{\mathbb{R}}
\newcommand{\A}{\mathscr{A}}
\newcommand{\B}{\mathscr{B}}
\newcommand{\C}{\mathscr{C}}
\newcommand{\T}{\mathscr{T}}
\newcommand{\E}{\mathscr{E}}

\DeclareMathOperator{\dist}{dist}
\long\def\comment#1{}

\newtheorem{theorem}{Theorem}[section]
\newtheorem*{theorem*}{Theorem}
\def\Thm#1#2{\ifthenelse{\equal{#1}{*}}{\begin{theorem*}#2\end{theorem*}}
  {\begin{theorem}\label{T#1}#2\end{theorem}}}
\def\thm#1{Theorem~\ref{T#1}}

\newtheorem{proposition}[theorem]{Proposition}
\newtheorem*{proposition*}{Proposition}
\def\Prp#1#2{\ifthenelse{\equal{#1}{*}}{\begin{proposition*}#2\end{proposition*}}
             {\begin{proposition}\label{P#1}#2\end{proposition}}}
\def\prp#1{Proposition~\ref{P#1}}

\newtheorem{corollary}[theorem]{Corollary}
\newtheorem*{corollary*}{Corollary}
\def\Cor#1#2{\ifthenelse{\equal{#1}{*}}{\begin{corollary*}#2\end{corollary*}}
             {\begin{corollary}\label{C#1}#2\end{corollary}}}
\def\cor#1{Corollary~\ref{C#1}}

\newtheorem{lemma}[theorem]{Lemma}
\newtheorem*{lemma*}{Lemma}
\def\Lem#1#2{\ifthenelse{\equal{#1}{*}}{\begin{lemma*}#2\end{lemma*}}
             {\begin{lemma}\label{L#1}#2\end{lemma}}}
\def\lem#1{Lemma~\ref{L#1}}

\theoremstyle{definition}
\newtheorem{remark}[theorem]{Remark}
\newtheorem*{remark*}{Remark}
\def\Rem#1#2{\ifthenelse{\equal{#1}{*}}{\begin{remark*}\rm #2\end{remark*}}
             {\begin{remark}\label{R#1}\rm #2\end{remark}}}

\newcommand{\eq}[1]{\eqref{E#1}}
\def\Eq#1#2{\ifthenelse{\equal{#1}{*}}
  {\begin{equation*}\begin{aligned}[]#2\end{aligned}\end{equation*}}
  {\begin{equation}\begin{aligned}[]\label{E#1}#2\end{aligned}\end{equation}}}

\begin{document}

\title{Estimates for approximately Jensen convex functions}

\author{G\'abor M. Moln\'ar and Zsolt P\'ales}
\address{G\'abor M. Moln\'ar, Institute of Mathematics and Computer Sciences, University of Nyíregyháza, H-4400 Nyíregyháza, Sóstói út 31/B, Hungary}
\email{molnar.gabor@nye.hu}

\address{Zsolt P\'ales, Institute of Mathematics, University of Debrecen, H-4002 Debrecen, Pf.\ 400, Hungary}
\email{pales@science.unideb.hu}



\thanks{The research of the first author was supported by the Scientific Council of the University of Nyíregyháza and by the EKÖP-24 University Excellence Scholarship Program of the Ministry for Culture and Innovation from the source of the National Research, Development and Innovation Fund. \\\indent The research of the second author was supported by the K-134191 NKFIH Grant.}

\begin{abstract}
In this paper functions $f:D\to\mathbb{R}$ satisfying the inequality
$$
  f\Big(\frac{x+y}{2}\Big)\leq\frac12f(x)+\frac12f(y)
  +\varphi\Big(\frac{x-y}{2}\Big) \qquad(x,y\in D)
$$
are studied, where $D$ is a nonempty convex subset of a real linear space $X$ and $\varphi:\{\frac12(x-y) : x,y \in D\}\to\mathbb{R}$ is a so-called error function. In this situation $f$ is said to be \emph{$\varphi$-Jensen convex.} The main results show that for all $\varphi$-Jensen convex function $f:D\to\mathbb{R}$, for all rational $\lambda\in[0,1]$ and $x,y\in D$, the following inequality holds
$$
  f(\lambda x+(1-\lambda)y)
  \leq \lambda f(x)+(1-\lambda)f(y)+\sum_{k=0}^\infty
  \frac{1}{2^k}\varphi\big(\dist(2^k\lambda,\mathbb{Z})\cdot(x-y)\big).
$$
The infinite series on the right hand side is always convergent, moreover, for all rational $\lambda\in[0,1]$, it can be evaluated as a finite sum. 
\end{abstract}

\maketitle

\section{Introduction}

Throughout the paper let $\R$, $\R_+$, $\Q$, $\Z$ and $\N$ denote the sets of real, nonnegative real, rational, integer, and natural numbers, respectively. For brevity, we also denote the set of rational elements of a real interval $I$ by $I_\Q$, i.e., $I_\Q:=I\cap\Q$. In what follows, we always assume that $D$ is a nonempty convex subset of a real linear space $X$. We say that a function $f:D\to\R$ is \emph{Jensen convex on $D$} (or \emph{midconvex on $D$}) if 
\Eq{JC}{
  f\Big(\frac{x+y}{2}\Big)\leq\frac12f(x)+\frac12f(y)
   \qquad(x,y\in D).
}
From the theory of Jensen convex functions, we recall the following classical result which can be found in the monograph \cite[Theorem 5.3.5.]{Kuc85} of M.\ Kuczma (see also \cite{HarLitPol34}, \cite{NicPer18}).

\Thm{QC}{Let $f:D\to\R$ be Jensen convex (i.e., assume that \eq{AJC} holds with $\varphi\equiv0$. Then $f$ is $\Q$-convex, that is, for all $\lambda\in[0,1]_\Q$ and $x,y\in D$,
\Eq{QC}{
f(\lambda x + (1-\lambda)y)\leq \lambda f(x) + (1-\lambda)f(y).
}}

A counterpart of this result is the following theorem.

\Thm{QCC}{If $\lambda\in[0,1]\setminus\Q$ and $x,y$ are distinct elements of $D$, then there exists a Jensen convex function $f:D\to\R$ such that the reversed strict inequality is valid in \eq{QC}.}

\begin{proof}
Let $x,y$ be distinct elements of $D$, $\lambda\in[0,1]\setminus\Q$, and let $H$ be a Hamel base of $X$ over $\Q$ (cf.\ \cite{Ham05}), which contains the $\Q$-independent vectors $x-y$ and $\lambda (x-y)$. Let $A:X\to\R$ be an additive function such that $A(x-y)=0$ and $A(\lambda (x-y))=1$ and define $f:D\to\R$ by $f(u):=A(u-y)$. Then, $f$ is obviously a Jensen convex function, on the other hand,
\Eq{*}{
   1=A(\lambda (x-y))
   &=f(\lambda x + (1-\lambda)y)\\
   &> \lambda f(x) + (1-\lambda)f(y)
   =\lambda A(x-y)+(1-\lambda)A(y-y)=0,
}
which shows that \eq{QC} does not hold.
\end{proof}

The above result shows that the inequality \eq{QC} does not remain valid for all $\lambda\in[0,1]$ without any further assumption on $f$. The continuity of $f$ (along the lines, i.e., the 1-dimensional affine subsets of $D$) is a sufficient condition. A more powerful sufficient condition was discovered by Bernstein and Doetsch \cite{BerDoe15}.

\Thm{BD}{Let $X$ be a topological linear space and let $D\subseteq X$ be a nonempty open and convex set. Let $f:D\to\R$ be a Jensen convex function which is locally bounded from above at a point of $D$. Then $f$ is continuous and \eq{QC} holds for all $\lambda\in[0,1]$ and $x,y\in D$.}

The investigation of approximate convexity started with the seminal papers of \cite{Gre52}, \cite{HyeUla52}. Since then, there have been several papers published dealing with various generalizations of an approximate or a strong version of Jensen convexity as well as of convexity for real-valued and also for set-valued functions. See, for instance the papers \cite{Ada16,AzoGimNikSan11,Bor08,GerNik11,GilGonNikPal17,GonNikPalRoa14,Haz05a,Haz07b,Haz11,HazPal04,HazPal05,Jov96,Lac99,LeiMerNikSan13,MakPal10b,MakPal13b,MerNik10,Mro01,MurTabTab12,NgNik93,Nik89,Nik89b,NikPal11,Pal00b,Pal03a,Pol96,Raj15,RajWas11,TabTab09b,TabTab09a, TabTabZol10b,TabTabZol10a}.

To recall a notion of approximate Jensen convexity (introduced in the paper \cite{MakPal12a}), given a convex set $D\subseteq X$, we define the difference set $D_\Delta$ by $D_\Delta:= \{x-y : x,y \in D\}$. Then $D_\Delta$ is also a convex set which is symmetric with respect to the origin, i.e., it satisfies $D_\Delta=-D_\Delta$ and hence it contains the origin of $X$. Let $\varphi:\frac12D_\Delta\to\R$ be a given function, which we will refer to as an \emph{error function}. We say that a function $f:D\to\R$ is \emph{$\varphi$-Jensen convex on $D$} (or \emph{$\varphi$-midconvex on $D$}) if 
\Eq{AJC}{
  f\Big(\frac{x+y}{2}\Big)\leq\frac12f(x)+\frac12f(y)
  +\varphi\Big(\frac{x-y}{2}\Big) \qquad(x,y\in D).
}
Interchanging the variables $x$ and $y$ in \eq{AJC}, we get
\Eq{*}{
  f\Big(\frac{x+y}{2}\Big)\leq\frac12f(x)+\frac12f(y)
  +\varphi\Big(\frac{y-x}{2}\Big) \qquad(x,y\in D),
}
Combining this inequality with \eq{AJC}, it follows that
\Eq{*}{
  f\Big(\frac{x+y}{2}\Big)\leq\frac12f(x)+\frac12f(y)
  +\min\Big(\varphi\Big(\frac{x-y}{2}\Big),\varphi\Big(\frac{y-x}{2}\Big)\Big) \qquad(x,y\in D).
}
Therefore, we may assume that $\varphi$ is an even function, i.e., for all $u\in\frac12D_\Delta$, we have $\varphi(u)=\varphi(-u)$. If $x=y$, then \eq{AJC} is equivalent to the inequality $0\leq\varphi(0)$, thus we may also assume that $\varphi(0)=0$.

If $\varphi$ is nonnegative, then we usually speak about the \emph{approximate convexity of $f$}. On the other hand, if $\varphi$ is nonpositive, then the inequality \eq{AJC} implies the standard Jensen convexity of $f$, therefore, in this case we may speak about \emph{strong Jensen convexity of $f$}. A tipical choice in the setting of strong convexity is when $\varphi(u)=-\varepsilon\|u\|^2$ for some positive $\varepsilon$ and $X$ is an inner product space.  

In the case when $\varphi$ is nonnegative, the zero function $f=0$ is always $\varphi$-Jensen convex on $D$. However, if $\varphi$ can take negative values, then there may not exist a $\varphi$-Jensen convex function over $D$. Motivated by this, an even function $\varphi:\frac12D_\Delta\to\R$ which satisfies $\varphi(0)=0$ and for which there exists at least one $\varphi$-Jensen convex function over $D$ will be called an \emph{admissible error function}. In what follows, let $\A(D_\Delta)$ denote the family of all admissible error functions $\varphi:\frac12D_\Delta\to\R$.

The basic problem related to a $\varphi$-Jensen convex function $f:D\to\R$ is to deduce its further approximate convexity properties. More precisely, for a given $\varphi\in\A(D_\Delta)$, we search for an \emph{extended error function} $\Phi:[0,1]\times D_\Delta\to(-\infty,\infty]$ such that the inequality 
\Eq{AC}{
  f(\lambda x+(1-\lambda)y)\leq \lambda f(x)+(1-\lambda)f(y)+\Phi(\lambda,x-y)
  \qquad (\lambda\in[0,1],\,x,y\in D).
}
be valid for all $\varphi$-Jensen convex functions $f:D\to\R$.
If this inequality holds, then we say that $f$ is a \emph{$\Phi$-convex function on $D$}. Clearly, the extended error function $\Phi:[0,1]\times D_\Delta\to(-\infty,\infty]$ defined by
\Eq{*}{
  \Phi(\lambda,u)
  :=\begin{cases}
    \varphi\big(\tfrac12u\big) &\mbox{if } \lambda=\tfrac12 \mbox{ and }u\in D_\Delta,\\[2mm]
    0 &\mbox{if } \lambda\in\{0,1\} \mbox{ and } u\in D_\Delta 
    \mbox{ or } \lambda\in [0,1] \mbox{ and } u=0,\\[2mm]
    \infty &\mbox{otherwise}
    \end{cases}
}
is a trivial solution to this problem, because then \eq {AC} is equivalent to \eq{AJC}. It is not obvious if there is an extended error function satisfying \eq{AC} which is finite for $\lambda\in\big\{\tfrac13,\tfrac23\big\}$. To see that this is possible, let $x,y\in D$ be fixed. Now observe that $\tfrac{2x+y}{3}$ is the midpoint of the segment $\big[x,\frac{x+2y}{3}\big]$ and $\tfrac{x+2y}{3}$ is the midpoint of the segment $\big[\frac{2x+y}{3},y\big]$. Therefore, applying \eq{AJC} twice in these settings, we can get that
\Eq{*}{
   f\Big(\frac{2x+y}{3}\Big)&\leq \frac12f(x)+\frac12f\Big(\frac{x+2y}{3}\Big) 
   +\varphi\Big(\frac{x-y}{3}\Big),\\[2mm]
   f\Big(\frac{x+2y}{3}\Big)&\leq \frac12f\Big(\frac{2x+y}{3}\Big)+\frac12f(y) 
   +\varphi\Big(\frac{x-y}{3}\Big).
}
Using the second inequality, the first inequality yields that
\Eq{*}{
  f\Big(\frac{2x+y}{3}\Big)&\leq \frac12f(x)+\frac12\bigg(\frac12f\Big(\frac{2x+y}{3}\Big)+\frac12f(y)+\varphi\Big(\frac{x-y}{3}\Big)\bigg) 
   +\varphi\Big(\frac{x-y}{3}\Big). 
}
After simplifications, we get
\Eq{*}{
  f\Big(\frac{2x+y}{3}\Big)
  &\leq \frac23f(x)+\frac13f(y) 
   +2\varphi\Big(\frac{x-y}{3}\Big). 
}
In a similar manner, we can also obtain that
\Eq{*}{
  f\Big(\frac{x+2y}{3}\Big)
  &\leq \frac13f(x)+\frac23f(y) 
   +2\varphi\Big(\frac{x-y}{3}\Big). 
}
Thus, we have proved that \eq{AC} holds with 
the extended error function $\Phi:[0,1]\times D_\Delta\to(-\infty,\infty]$ defined by
\Eq{*}{
  \Phi(\lambda,u)
  :=\begin{cases}
    \varphi\big(\tfrac12u\big) &\mbox{if } \lambda=\tfrac12 \mbox{ and }u\in D_\Delta,\\[2mm]
    2\varphi\big(\tfrac13u\big) &\mbox{if } \lambda\in\big\{\tfrac13,\tfrac23\big\} \mbox{ and }u\in D_\Delta,\\[2mm]
    0 &\mbox{if } \lambda\in\{0,1\} \mbox{ and } u\in D_\Delta 
    \mbox{ or } \lambda\in [0,1] \mbox{ and } u=0,\\[2mm]
    \infty &\mbox{otherwise}.
    \end{cases}
}

In the context of $\varphi$-Jensen convexity, under regularity assumptions that are similar to those in the theorem of Bernstein and Doetsch, the following result was established in \cite{MakPal12a}.

\Thm{MP}{Let $X$ be a topological linear space and let $D\subseteq X$ be a nonempty open and convex set. Let $f:D\to\R$ be locally bounded from above at a point of $D$ and let $\varphi:\tfrac12 D_\Delta\to\R_+$ be nondecreasing. Then $f$ is $\varphi$-Jensen convex on $D$ if and only if, for all $\lambda\in[0,1]$ and $x,y\in D$,
\Eq{*}{
f(\lambda x + (1-\lambda)y)\leq \lambda f(x) + (1-\lambda)f(y) + \sum_{n=0}^\infty
  \frac{1}{2^n}\varphi\big(d_\Z(2^n\lambda)(x-y)\big),
}
where the function $d_\Z:\R\to\R$ is defined by
\Eq{dZ}{
  d_\Z(x):=\dist(x,\Z):=\inf\{|x-k|\mid k\in\Z\} \qquad (x\in\R).
}}

The basic question concerning the above result is that what happens if we drop the local upper boundedness assumption on $f$ and also the monotonicity property of $\varphi$. These questions will be addressed in Section 3 of this paper. In Section 2, we collect some material which will be instrumental for establishing our results.

\section{Auxiliary results}

Motivated by \thm{MP}, we first establish some properties of the function $d_\Z$ defined by \eq{dZ} and of the map $\pmb\langle\cdot\pmb\rangle:\Q\to\N$ given by
\Eq{*}{
  \pmb\langle\lambda\pmb\rangle
  :=\min\{m\in\N\mid m\lambda\in\Z\}.
}
Observe that, for all $n\in\Z$, $x\in\R$, and $\lambda\in\Q$ the following obvious properties hold:
\Eq{*}{
  d_\Z(n\pm x)=d_\Z(x),
  \qquad
  \pmb\langle n\pm\lambda\pmb\rangle
  = \pmb\langle\lambda\pmb\rangle \qquad\mbox{and}\qquad
  \pmb\langle n\lambda\pmb\rangle
  \leq \pmb\langle\lambda\pmb\rangle.
}
For the function $d_\Z$, we have the following
easy-to-see properties
\Eq{ets}{
  0\leq d_\Z(x)\leq\tfrac12 \quad (x\in\R),\qquad
  d_\Z(x)=\min(x,1-x)\quad (x\in[0,1]).
}
The next lemma establishes a ``recursive'' formula for $d_\Z$.

\Lem{rec}{For all $x\in\R$, 
\Eq{d1}{
  d_\Z(2x)
  =d_\Z(2d_\Z(x))
  =\min\big(2d_\Z(x),1-2d_\Z(x)\big)
  =\begin{cases}
    2d_\Z(x) &\mbox{if } d_\Z(x)\leq\frac{1}{4},\\[2mm]
    1-2d_\Z(x) &\mbox{if } d_\Z(x)>\frac{1}{4}.
   \end{cases}
}
Furthermore, for all $k\in\N$ and $x\in\R$,
\Eq{dk}{
  d_\Z(2^kx)
  =\big(d_\Z\circ (2d_\Z)^k\big)(x).
}
}

\begin{proof} We need to show the first equality in \eq{d1}, the remaining two equalities easily follow from \eq{ets}. Let $x\in\R$ be fixed. To prove the first equality in \eq{d1}, we have to distinguish two cases.

If $d_\Z(x)\leq\frac{1}{4}$, then there exists an $n\in\N$ such that $n-\frac{1}{4}\leq x\leq n+\frac{1}{4}$. This implies that $2n-\frac{1}{2}\leq 2x\leq 2n+\frac{1}{2}$, therefore
\Eq{*}{
  d_\Z(2x)=|2x-2n|=2|x-n|=2d_\Z(x).
}

If $d_\Z(x)>\frac{1}{4}$, then there exists an $n\in\N$ such that $n+\frac{1}{4}<x<n+\frac{3}{4}$. Therefore, we have that $(2n+1)-\frac{1}{2}< 2x< (2n+1)+\frac{1}{2}$. Consequently,
\Eq{*}{
  d_\Z(2x)=|2x-(2n+1)|
  =1-2\big(\tfrac12-|n+\tfrac12-x|\big)=1-2d_\Z(x).
}
The equality \eq{dk} can be obtained by induction with respect to $k$. The equality holds for $k=1$ according to the first part of this lemma. If it holds for some $k$, then applying \eq{d1} first and then the inductive hypothesis, we get
\Eq{*}{
  d_\Z(2^{k+1}x)&=d_\Z(2\cdot 2^{k}x)
  =d_\Z\big(2d_\Z(2^kx)\big)\\
  &=d_\Z\big(2\big(d_\Z\circ (2d_\Z)^k\big)(x)\big)
  =d_\Z\big((2d_\Z)^{k+1}(x)\big)
  =\big(d_\Z\circ (2d_\Z)^{k+1}\big)(x).
}
This completes the proof of the lemma. 
\end{proof}

The following result and its subsequent corollary will play an important role in the comparison of the various upper bounds obtained for approximate convexity with rational convex combinations.

\Lem{psi}{Let $\psi:\big[0,\frac12\big]\to\R$ be bounded [continuous] and define the function $T_\psi:[0,1]\to\R$ by
\Eq{Tpsi}{
  T_\psi(\lambda):=\sum_{k=0}^\infty \frac{1}{2^k}
  \psi\big(d_\Z(2^k\lambda)\big).
}
Then $f=T_\psi$ is the unique bounded [continuous] solution of the functional equation
\Eq{TE}{
  f(\lambda)
  =\begin{cases}
    \frac12 f(2\lambda)+\psi(\lambda) &\mbox{if }\lambda\in\big[0,\frac12\big],
    \\[3mm]
    \frac12 f(2\lambda-1)+\psi(1-\lambda) &\mbox{if }\lambda\in\big(\frac12,1\big]
   \end{cases}
}
also satisfying the condition $f(0)=f(1)$.}

\begin{proof} The sequence of functions $\lambda\mapsto \psi\big(d_\Z(2^k\lambda)\big)$ is [continuous] and is uniformly bounded by $\|\psi\|_\infty$ (over $\big[0,\frac12\big]$). Therefore, the series on the right hand side of \eq{Tpsi} is majored by the numerical series $\sum_{k=0}^\infty \frac{1}{2^k}\|\psi\|_\infty=2\|\psi\|_\infty<\infty$. This yields that the series is uniformly convergent, whence it follows that the sum of the series is a bounded [continuous] function, i.e., $T_\psi$ is bounded [continuous]. Furthermore, we also have that
\Eq{*}{
  T_\psi(0)=\sum_{k=0}^\infty \frac{1}{2^k}
  \psi\big(d_\Z(2^k\cdot 0)\big)
  =\sum_{k=0}^\infty \frac{1}{2^k}
  \psi\big(0\big)
  =\sum_{k=0}^\infty \frac{1}{2^k}
  \psi\big(d_\Z(2^k\cdot 1)\big)=T_\psi(1).
}
If $\lambda\in\big[0,\frac12\big]$, then
\Eq{*}{
  T_\psi(2\lambda)=\sum_{k=0}^\infty \frac{1}{2^k}
  \psi\big(d_\Z(2^{k+1}\lambda)\big)
  =2\sum_{k=1}^\infty \frac{1}{2^{k}}
  \psi\big(d_\Z(2^{k}\lambda\big)\big)
  =2\big(T_\psi(\lambda)-\psi(\lambda)\big),
}
while, for $\lambda\in\big(\frac12,1\big]$,
\Eq{*}{
  T_\psi(2\lambda-1)&=\sum_{k=0}^\infty \frac{1}{2^k}
  \psi\big(d_\Z(2^{k+1}\lambda-2^k)\big)\\
  &=2\sum_{k=1}^\infty \frac{1}{2^{k}}
  \psi\big(d_\Z(2^{k}\lambda)\big)
  =2\big(T_p(\lambda)-\psi(1-\lambda)\big).
}
These two equalities show that the function $f=T_\psi$ satisfies the functional equation \eq{TE}.

We now show that $f=T_\psi$ is the unique solution of the functional equation \eq{TE} satisfying also the condition $f(0)=f(1)$. 

Let $\B\big([0,1],\R\big)$ and 
$\C\big([0,1],\R\big)$ denote normed spaces of bounded and continuous real valued functions defined on the interval $[0,1]$ (equipped with the supremum norm), respectively. Then these spaces are Banach spaces, furthermore, the functions belonging to theses spaces and satisfying the condition $f(0)=f(1)$ form closed subspaces, which we will denote by $\B_0\big([0,1],\R\big)$ and 
$\C_0\big([0,1],\R\big)$, respectively.

We point out that the right hand side of the functional equation \eq{TE} defines a contraction on these Banach subspaces. To see this, for a bounded [continuous] function $f:[0,1]\to\R$, define the transformed function $\T_\psi(f):[0,1]\to\R$ by
\Eq{*}{
  \T_\psi(f)(\lambda)
  =\begin{cases}
    \frac12 f(2\lambda)+\psi(\lambda) &\mbox{if }\lambda\in\big[0,\frac12\big],
    \\[3mm]
    \frac12 f(2\lambda-1)+\psi(1-\lambda) &\mbox{if }\lambda\in\big(\frac12,1\big].
   \end{cases}
}
It is easy to see that if $\psi$ is bounded [continous], then $\T_\psi$ maps the space $\B_0\big([0,1],\R\big)$ [resp. $\C_0\big([0,1],\R\big)$] into itself. Indeed, if $f(0)=f(1)$, then
\Eq{*}{
  \T_\psi(f)(0)
  =\frac12 f(0)+\psi(0)
  =\frac12 f(1)+\psi(0)
  =\T_\psi(f)(1).
}

If $\psi$ and $f$ are bounded, then $\T_\psi(f)$ is also bounded by $\frac12 \|f\|_\infty+\|\psi\|_\infty$, which shows that $\T_\psi(f)$ maps the subspace $\B_0\big([0,1],\R\big)$ into itself.

If $\psi$ and $f$ are continuous, then the continuoity of $\T_\psi(f)$ on $\big[0,\frac12\big)\cup\big(\frac12,1\big]$ is obvious. It remains to check the continuity at $\lambda=\frac12$.
\Eq{*}{
  \lim_{\lambda\to\frac12-}&\T_\psi(f)(\lambda)
  =\lim_{\lambda\to\frac12-}\Big(\tfrac12 f(2\lambda)+\psi(\lambda)\Big)=\tfrac12f(1)+\psi\big(\tfrac12\big)\\
  &=\tfrac12f(0)+\psi\big(\tfrac12\big)
  =\lim_{\lambda\to\frac12+}\Big(\tfrac12 f(2\lambda-1)+\psi(1-\lambda)\Big)
  =\lim_{\lambda\to\frac12+}\T_\psi(f)(\lambda).
}

Finally, we show that $\T_\psi$ is a contraction with contraction factor $q=\frac12$.
To see this, let $f,g\in\B_0\big([0,1],\R\big)$.
Then, for $\lambda\in\big[0,\frac12\big]$,
\Eq{*}{
  \big|\T_\psi(f)(\lambda)-\T_\psi(g)(\lambda)\big|
  =\big|\tfrac12 f(2\lambda)-\tfrac12 g(2\lambda)\big|
  \leq\tfrac12\|f-g\|_\infty,
}
while, for $\lambda\in\big(\frac12,1\big]$,
\Eq{*}{
  \big|\T_\psi(f)(\lambda)-\T_\psi(g)(\lambda)\big|
  =\big|\tfrac12 f(2\lambda-1)-\tfrac12 g(2\lambda-1)\big|
  \leq\tfrac12\|f-g\|_\infty.
}
These two inequalities yield that
\Eq{*}{
  \big\|\T_\psi(f)-\T_\psi(g)\big\|_\infty
  \leq\tfrac12\|f-g\|_\infty,
}
which shows that $\T_\psi$ is a $\frac12$-contraction on $\B_0\big([0,1],\R\big)$, indeed.

It is obvious that $f:[0,1]\to\R$ is a solution of the functional equation \eq{TE} if and only if $f$ is a fixed point of $\T_\psi$. According to the Banach Fixed Point Theorem, $\T_\psi$ has a unique fixed point, and this completes the proof of the unique solvability of the functional equation \eq{TE} in the subspaces $\B_0\big([0,1],\R\big)$ [and $\C_0\big([0,1],\R\big)$].
\end{proof}

The function $T_\psi$ defined in \eq{Tpsi} is called the \emph{Takagi-type function related to $\psi$}. In the particular case when $\psi(x)=x$ it is the classical continuous but nowhere differentiable function introduced by Takagi \cite{Tak03}. In what follows, we compute $T_\psi$ for the case when $\psi(x)=x^2$. We note that this is the only power function for which $T_\psi$ is smooth on $\R\setminus\Z$.

\Cor{psi}{For all $\lambda\in[0,1]$, 
\Eq{p=2}{
  \sum_{k=0}^\infty \frac{1}{2^k}
  d_\Z(2^k\lambda)^2
  =\lambda(1-\lambda).
}}

\begin{proof} We apply the \lem{psi} for the function $\psi(t):=t^2$ \ ($t\in\big[0,\frac12\big]$). According to this lemma, $f=T_\psi$ is the unique solution continous of the functional equation
\Eq{fe}{
  f(\lambda)
  =\begin{cases}
    \frac12 f(2\lambda)+\lambda^2 &\mbox{if }\lambda\in\big[0,\frac12\big],
    \\[3mm]
    \frac12 f(2\lambda-1)+(1-\lambda)^2 &\mbox{if }\lambda\in\big(\frac12,1\big],
   \end{cases}
}
which also satisfies $f(0)=f(1)$. To verify the equality \eq{p=2}, it is sufficient to show that the function $f(\lambda)=\lambda(1-\lambda)$ is also a solution to \eq{fe}. 

If $\lambda\in\big[0,\frac12\big]$, then
\Eq{*}{
  f(\lambda)=\lambda(1-\lambda)
  =\tfrac12 (2\lambda)(1-2\lambda)+\lambda^2
  =\tfrac12 f(2\lambda)+\lambda^2,
}
while, for $\lambda\in\big(\frac12,1\big]$,
\Eq{*}{
  f(\lambda)=\lambda(1-\lambda)
  =\frac12 (2\lambda-1)(1-(2\lambda-1))+(1-\lambda)^2
  =\frac12 f(2\lambda-1)+(1-\lambda)^2.
}
These two identities show that the function $f(\lambda)=\lambda(1-\lambda)$ also satisfies the functional equation \eq{fe}. By the unique solvability of \eq{fe}, it follows that the equality \eq{p=2} holds for all $\lambda\in[0,1]$.
\end{proof}

In what follows, let the symbol $\phi$ denote Euler's totient function. Recall that, according to Euler's celebrated result,
$n\mid (a^{\phi(n)}-1)$ if $n$ and $a$ are coprime natural numbers. The next statement sharpens Euler's theorem in some particular cases. The authors are grateful to Professor Lajos Hajdu for providing some ideas for its proof.

\Lem{ET}{Let $a,n\in\N$ be coprime numbers such that $n\geq3$ is an odd number, and let $\ell:=\frac12\phi(n)$. Then either $n\mid (a^\ell-1)$ or $n\mid(a^\ell+1)$.}

\begin{proof} We note that the assumption $n\geq3$ yields that $\phi(n)$ is even, i.e., $\ell$ is an integer. 

Assume first that $n$ is a prime power, say $n=p^k$, where $p$ is an odd prime and $k>1$ is an integer. According to Euler's Theorem, $n\mid(a^{\phi(n)}-1)$, i.e., $n\mid (a^\ell-1)(a^\ell+1)$. The greatest common divisor of the numbers $(a^\ell-1)$ and $(a^\ell+1)$ is $2$ if both of them are even, otherwise they are coprime because their difference is divisible by $2$. 

Since $p$ is an odd prime, then $p$ cannot divide both of $(a^\ell-1)$ and $(a^\ell+1)$, therefore $n=p^k$ must be a divisor either of $(a^\ell-1)$ or of $(a^\ell+1)$.

Assume now that $n$ is not a power of a prime number. Then $n$ is of the form $n=km$, where $k,m\geq3$ and $k,m$ are coprime odd numbers which are also coprime to $a$. To complete the proof of the lemma, we show that $k$ and $m$ are divisors of $(a^\ell-1)$. By the multiplicativity of Euler's totient function, we have that $\phi(n)=\phi(k)\phi(m)$. 
Therefore, $\ell=\phi(k)\cdot\frac{\phi(m)}{2}=\phi(m)\cdot\frac{\phi(k)}{2}$. Thus, $\phi(k)$ and $\phi(m)$ are divisors of $\ell$. This implies that $a^{\phi(k)}-1$ and $a^{\phi(m)}-1$ are divisors of $a^\ell-1$. On the other hand, acording to Euler's theorem, $k\mid (a^{\phi(k)}-1)$ and $m\mid (a^{\phi(m)}-1)$, which then yields that $k,m\mid (a^\ell-1)$.
Using that $k$ and $m$ are coprime numbers, we can conclude that $n=km\mid (a^\ell-1)$.
\end{proof}

The role of the function $\mu_n$ introduced in the following lemma will be clarified in \prp{mk} below.

\Lem{mun}{Let $n\geq3$ be an odd integer and let $M_n$ denote the set of those integers $m\in\big\{1,\dots,\frac{n-1}{2}\big\}$ that are coprime to $n$. Define the mapping $\mu_n:M_n\to\Z$ defined by 
\Eq{*}{
 \mu_n(m):=\min(2m,n-2m).
}
Then $\mu_n$ is a bijection of the set $M_n$, furthermore, for all $m\in M_n$ and $k\in\N$,
\Eq{munk}{
  \{\mu_n^k(m)\}
  = \Big(\big\{2^km-in\mid i\in\{0,\dots,2^{k-1}-1\}\big\}
  \cup \big\{in-2^km\mid i\in\{1,\dots,2^{k-1}\}\big\}\Big)\cap M_n,
}
and the equality $m=\mu_n^\ell(m)$ holds for all $m\in M_n$, where $\ell:=\frac12\phi(n)$.}

\begin{proof}
If $m\in M_n$, then $m$ is coprime to $n$, therefore $2m$ and $n-2m$ are also coprimes to $n$. On the other hand, the inequalities $1\leq m\leq\frac{n-1}{2}$ imply that $1\leq\min(2m,n-2m)=\mu_n(m)$ and $\mu_n(m)=\min(2m,n-2m)<\frac{2m+(n-2m)}{2}=\frac{n}{2}$ and hence $\mu_n(m)\leq\frac{n-1}{2}$. Thus, we have proved that $\mu_n(m)\in M_n$. To see that $\mu_n$ is injective, assume that, for some $m,m'\in M_n$, we have that $\mu_n(m)=\mu_n(m')$, i.e., $\min(2m,n-2m)=\min(2m',n-2m')$. This equality implies that one of the following equalities must be valid:
\Eq{*}{
 2m=2m',\qquad
 2m=n-2m',\qquad
 n-2m=2m',\qquad
 n-2m=n-2m'.
}
The second and third equalities contradict the oddness of $n$. Therefore, either the first or the last equality must hold, which imply that $m=m'$. This shows that $\mu_n$ is injective.

To prove the equality \eq{munk}, we first verify by induction on $k$ that
\Eq{munkk}{
  \mu_n^k(m)
  \in\big\{2^km-in\mid i\in\{0,\dots,2^{k-1}-1\}\big\}
  \cup\big\{in-2^km\mid i\in\{1,\dots,2^{k-1}\}\big\}.
}
For $k=1$, we have 
\Eq{*}{
  \mu_n(m)=\min(2m,n-2m)\in\{2m,n-2m\}
  =\{2m-0\cdot n\}\cup \{1\cdot n-2m\},
}
which is the inclusion \eq{munkk} in the particular case $k=1$.

Assume that \eq{munkk} is valid for some $k\in\N$. Then
\Eq{*}{
  &\mu_n^{k+1}(m)=\mu_n\big(\mu_n^k(m)\big)
  =\min\big(2\mu_n^k(m),n-2\mu_n^k(m)\big)\in\big\{2\mu_n^k(m),n-2\mu_n^k(m)\big\}\\
  &\subseteq \big\{2^{k+1}m-2in\mid i\in\{0,\dots,2^{k-1}-1\}\big\}
  \cup \big\{2in-2^{k+1} m\mid i\in\{1,\dots,2^{k-1}\}\big\}\\
     &\quad\cup\big\{(2i+1)n-2^{k+1}m\mid i\in\{0,\dots,2^{k-1}-1\}\big\}
  \cup \big\{2^{k+1}m-(2i-1)n\mid i\in\{1,\dots,2^{k-1}\}\big\}\\
  &=\big\{2^{k+1}m-in\mid i\in\{0,\dots,2^{k}-1\}\big\}
  \cup \big\{in-2^{k+1}m\mid i\in\{1,\dots,2^{k}\}\big\},
}
which proves the assertion for $k+1$. 

Since, $\mu_n$ maps $M_n$ into itself, we have that $\mu_n^k(m)\in M_n$ for all $k\in\N$ and $m\in M_n$, whence it follows that
\Eq{*}{
  \{\mu_n^k(m)\}
  \subseteq \Big(\big\{2^km-in\mid i\in\{0,\dots,2^{k-1}-1\}\big\}
  \cup \big\{in-2^km\mid i\in\{1,\dots,2^{k-1}\}\big\}\Big)\cap M_n.
}
To show that here the equality holds, it suffices to prove that the set on the right hand side of this inclusion is a singleton.
Let $p$ be an element of the right hand side of the above inclusion. Then either $p=2^km-in$ for some $i\in\{0,\dots,2^{k-1}-1\}$ or $p=in-2^km$ for some $i\in\{1,\dots,2^{k-1}\}$.
In both cases, due to the oddness of $n$, it follows that $p$ and $n$ are coprime. Therefore $p$ belongs to $M_n$ if and only if the inequalities $0<p<\frac{n}{2}$ holds. In the first, resp.\ in the second case these inequalities hold if and only if 
\Eq{ii}{
2^k\frac{m}{n}-\frac{1}{2}<i<2^k\frac{m}{n}, \qquad\mbox{resp.}\qquad 2^k\frac{m}{n}<i<2^k\frac{m}{n}+\frac{1}{2}.
}
Thus, the number $i$ has to be the unique integer in the open interval $\big]2^k\frac{m}{n}-\frac{1}{2},2^k\frac{m}{n}+\frac{1}{2}\big[$ (whose length is equal to 1). In the case when the first inequality is valid in \eq{ii}, then $p$ is of the form $p=2^km-in$, in the other case $p=in-2^km$. Thus $p$ is uniquely determined and this verifies the equality \eq{munk}.

Observe that (by the definition of the totient function) the cardinality of the set $M_n$ equals $\ell=\frac12\phi(n)$. For every $m\in M_n$, we have that $m,\mu_n(m),\mu_n^2(m),\dots,\mu_n^\ell(m)\in M_n$, therefore, according to the pigeonhole principle, there exist $0\leq i<j\leq\ell$ such that $\mu_n^i(m)=\mu_n^j(m)$. Consequently, $m=\mu_n^{j-i}(m)$, i.e., $m=\mu_n^{k}(m)$ for some $k\in\{1,\dots,\ell\}$. We may assume that $k$ is the smallest such exponent. 
Using the inclusion \eq{munk}, it follows that if $m=\mu_n^k(m)$, then one of the following equalities must hold:
\Eq{*}{
  (2^k-1)m=in \qquad(i\in\{1,\dots,2^{k-1}-1\}),\qquad
  (2^k+1)m=in \qquad(i\in\{1,\dots,2^{k-1}\}).
}
\comment{i.e.,
\Eq{*}{
  \frac{m}{n}\in \bigg\{\frac{1}{2^k-1},\dots,\frac{2^{k-1}-1}{2^k-1}\bigg\}
  \cup \bigg\{\frac{1}{2^k+1},\dots,\frac{2^{k-1}}{2^k+1}\bigg\}.
}}
Using that $m$ is coprime to $n$, it follows that $m\mid i$ and either $n\mid (2^k-1)$ or $n\mid (2^k+1)$.

In what follows, we are going to prove that $k\mid\ell$. If this were not be true, then $\ell=dk+r$ would hold for some $d,k\in\N$ with $0<r<k$. According to \lem{ET}, we have that either $n\mid (2^\ell-1)$ or $n\mid (2^\ell+1)$. Thus, we can distinguish the following four cases:

Case I: $n\mid (2^k-1)$ and $n\mid (2^\ell-1)$. Then $n\mid(2^{dk}-1)$. Thus $n\mid(2^\ell-2^{dk})=2^{dk}(2^r-1)$, hence $n\mid(2^r-1)$ (because $n$ is odd).

Case II: $n\mid (2^k-1)$ and $n\mid (2^\ell+1)$. Then $n\mid(2^{dk}-1)$. Thus $n\mid(2^\ell+2^{dk}\big)=2^{dk}(2^r+1)$, hence $n\mid(2^r+1)$.

Case III: $n\mid (2^k+1)$ and $n\mid (2^\ell-1)$. Then $n\mid(2^{dk}-(-1)^d)$. Thus $n\mid(2^\ell+(-1)^d2^{dk})=2^{dk}(2^r-(-1)^d)$, hence $n|(2^r+(-1)^d)$.

Case IV: $n\mid (2^k+1)$ and $n\mid (2^\ell+1)$. Then $n\mid(2^{dk}-(-1)^d)$. Thus $n\mid(2^\ell-(-1)^d2^{dk})=2^{dk}(2^r+(-1)^d)$, hence $n\mid(2^r-(-1)^d)$.

These cases show that either $n\mid (2^r-1)$ or $n\mid (2^r+1)$.
This implies that
\Eq{*}{
  \frac{m}{n}\in \bigg\{\frac{1}{2^r-1},\dots,\frac{2^{r-1}-1}{2^r-1}\bigg\}
  \cup \bigg\{\frac{1}{2^r+1},\dots,\frac{2^{r-1}}{2^r+1}\bigg\}.
}
The above inclusion yield that
\Eq{*}{
  m\in \Big(\big\{2^rm-in\mid i\in\{1,\dots,2^{r-1}-1\}\big\}
  \cup \big\{in-2^rm\mid i\in\{1,\dots,2^{r-1}\}\big\}\Big)\cap M_n,
}
which, according to the second assertion of this lemma, implies that $m=\mu_n^r(m)$. Since $r<k$, this contradicts the minimality of $k$. The contradiction so obtained shows that $k$ is a divisor of $\ell$ and hence $m=\mu_n^\ell(m)$. (In other words, $\mu_n^\ell$ is the identical self map of the set $M_n$).
\end{proof}

In the next statement, we establish a connection between the functions $d_\Z$ and $\mu_n$. 

\Prp{mk}{Let $n\geq 3$ and let $M_n$ and $\mu_n$ denote the set and function introduced in \lem{mun}. Then, for $m\in M_n$ and $k\geq0$,
\Eq{mk}{
  d_\Z(2^{k}\tfrac{m}{n})=\tfrac{1}{n}\mu_n^k(m).
}}

\begin{proof}
We prove the statement of the proposition applying induction on $k$.

The equality \eq{mk} is obvious for $k=0$ because both sides are equal to $\tfrac{m}{n}$. For $k=1$, we can use formula \eq{d1} of \lem{rec} with $x=\tfrac{m}{n}$. Then
\Eq{*}{
  d_\Z(2\tfrac{m}{n})
  &=\min(2d_\Z(\tfrac{m}{n}),1-2d_\Z(\tfrac{m}{n}))\\
  &=\min(2\tfrac{m}{n},1-2\tfrac{m}{n})
  =\tfrac{1}{n}\min(2m,n-2m)=\tfrac{1}{n}\mu_n(m),
}
which yields that \eq{mk} is valid for $k=1$.

Assume now that the equality \eq{mk} holds for some $k\geq1$. Then, applying formula \eq{dk} with $x=\tfrac{m}{n}$ twice and the inductive hypothesis (also twice), we get
\Eq{*}{
  d_\Z(2^{k+1}\tfrac{m}{n})
  &=(d_\Z\circ(2d_\Z)^{k+1})(\tfrac{m}{n})
  =(d_\Z\circ(2d_Z)\circ(2d_\Z)^{k})(\tfrac{m}{n})\\
  &=d_\Z\big(2(d_Z\circ(2d_\Z)^{k})(\tfrac{m}{n})\big)
  =d_\Z\big(2d_\Z(2^{k}\tfrac{m}{n})\big)\\
  &=d_\Z\big(\tfrac{2}{n}\mu_n^k(m)\big)
  =\tfrac{1}{n}\mu_n(\mu_n^k(m))
  =\tfrac{1}{n}\mu_n^{k+1}(m).
}
This shows that the equality \eq{mk} holds for $k+1$ (instead of $k$), and this completes the induction.
\end{proof}

Finally, we prove a certain periodicity property of the sequence $d_\Z(2^k\lambda)$.

\Prp{kl}{Let $n\geq 3$ be an odd number, $m\in\Z$ be coprime to $n$ and let $\ell:=\frac12\phi(n)$. Then, for all nonnegative $k\in\Z$, 
\Eq{kl}{
  d_\Z(\tfrac{m}{n})
  =d_\Z(2^{k\ell}\cdot\tfrac{m}{n}).
}
Consequently, the sequence $\big(d_\Z(2^{k}\cdot\tfrac{m}{n})\big)_{k\geq0}$ is $\ell$-periodic.}

\begin{proof}
Then $\ell$ is an integer since $\phi(n)$ is even for $n\geq2$. The value of the function $d_\Z$ at $\frac{m}{n}$ is rational number of the form $\frac{m'}{n}\in\big(0,\frac12\big)$ such that either $n\mid (m-m')$ or $n\mid (m+m')$. Therefore, there exists $d\in\Z$ such that either $m=dn+m'$ or $m=dn-m'$. Since $m$ is coprime to $n$, it follows that $m'$ is also coprime to $m$. Therefore, $m'\in M_n$, and according to the last assertion of \lem{mun}, for all $k\geq0$, we have that $m'=\mu_n^{k\ell}(m')$.

In the case when $m=dn+m'$, using \prp{mk}, we get
\Eq{*}{
  d_\Z(\tfrac{m}{n})
  &=d_\Z(\tfrac{dn+m'}{n})
  =d_\Z(d+\tfrac{m'}{n})
  =d_\Z(\tfrac{m'}{n})
  =\tfrac{1}{n}m'
  =\tfrac{1}{n}\mu_n^{k\ell}(m')\\
  &=d_\Z(2^{k\ell}\cdot\tfrac{m'}{n})
  =d_\Z(2^{k\ell}\cdot\tfrac{m-dn}{n})
  =d_\Z(2^{k\ell}\cdot\tfrac{m}{n}-2^{k+\ell}d)
  =d_\Z(2^{k\ell}\cdot\tfrac{m}{n}).
}
In the case when $m=dn-m'$, using \prp{mk}, we get
\Eq{*}{
  d_\Z(\tfrac{m}{n})
  &=d_\Z(\tfrac{dn-m'}{n})
  =d_\Z(d-\tfrac{m'}{n})
  =d_\Z(-\tfrac{m'}{n})
  =d_\Z(\tfrac{m'}{n})
  =\tfrac{1}{n}m'
  =\tfrac{1}{n}\mu_n^{k\ell}(m')\\
  &=d_\Z(2^{k\ell}\cdot\tfrac{m'}{n})
  =d_\Z(2^{k\ell}\cdot\tfrac{dn-m}{n})
  =d_\Z(2^{k\ell}d-2^{k\ell}\cdot\tfrac{m}{n})
  =d_\Z(-2^{k\ell}\cdot\tfrac{m}{n})
  =d_\Z(2^{k\ell}\cdot\tfrac{m}{n}).
}
Thus, we have proved that \eq{kl} holds in both cases.
\end{proof}

\section{Main results}

In our first result, we show that $\varphi$-Jensen convexity implies further convexity properties for rational convex combinations of the variables. 

\Thm{AQC}{Let $\varphi\in\A(D_\Delta)$. If $f:D\to\R$ is $\varphi$-Jensen convex, then, for all $n,k\in\N$ and $x,y\in D$,
\Eq{nk}{
  f\Big(\dfrac{nx+ky}{n+k}\Big)
  \leq \frac{n}{n+k}f(x)+\frac{k}{n+k}f(y)+
  kn\varphi\bigg(\frac{x-y}{n+k}\bigg).
}}

\begin{proof} We are going to show this inequality by induction with respect to $m=n+k$. 

If $m=n+k=2$, then $n=k=1$ and \eq{nk} is equivalent to the $\varphi$-Jensen convexity of $f$.

Let $m\geq2$ and assume that the inequality \eq{nk} is valid for all $n,k\in\N$ such that $n+k\leq m$. Let $n,k\in\N$ such that $n+k=m+1$. Then
\Eq{*}{
	\max(n,k)\geq\frac{m+1}{2}\geq\frac{3}{2},
}
so $\max(n,k)\geq2$. Without loss of generality, we may assume  that $n\geq2$. Observe that
\Eq{*}{
	\frac{nx+ky}{n+k}=\frac{m-k}{m}\cdot x+\frac{k}{m}\cdot\frac{x+my}{n+k}.
}
Therefore, by the inductive hypothesis, we can apply the inequality \eq{nk} for $(m-k,k)\in\N^2$ and $\big(x,\frac{x+my}{m+1}\big)\in D^2$, and we obtain
\Eq{nxky}{
	f\bigg(\frac{nx+ky}{n+k}\bigg)
	&\leq\frac{m-k}{m}f(x)+\frac{k}{m}f\bigg(\frac{x+my}{n+k}\bigg)+(m-k)k\varphi\bigg(\frac{1}{(m-k)+k}\Big(x-\frac{x+my}{n+k}\Big)\bigg)\\
	&\leq\frac{m-k}{m}f(x)+\frac{k}{m}f\bigg(\frac{x+my}{n+k}\bigg)+(m-k)k\varphi\bigg(\frac{x-y}{n+k}\bigg).
}
On the other hand, we also have that
\Eq{*}{
	\frac{x+my}{n+k}
	=\frac{1}{n}\cdot\frac{nx+ky}{n+k}+\frac{n-1}{n}\cdot y.
}
Notice that $1+(n-1)=n=m+1-k\leq m$.
Thus, by the inductive assumption, we can apply the inequality \eq{nk} for $(1,n-1)\in\N^2$ and $\big(\frac{nx+ky}{n+k},y\big)\in D^2$, and we get
\Eq{*}{
	f\bigg(\frac{x+my}{n+k}\bigg)&\leq\frac{1}{n}f\bigg(\frac{nx+ky}{n+k}\bigg)+\frac{n-1}{n}f(y)+1\cdot (n-1)\varphi\bigg( \frac{1}{1+(n-1)}\Big(\frac{nx+ky}{n+k}-y \Big)\bigg)\\
	&=\frac{1}{n}f\bigg(\frac{nx+ky}{n+k}\bigg)+\frac{n-1}{n}f(y)+(n-1)\varphi\bigg(\frac{x-y}{n+k}\bigg).
}
Using this inequality, \eq{nxky} yields that
\Eq{*}{
	f\bigg(\frac{nx+ky}{n+k}\bigg)
	\leq\frac{m-k}{m}f(x)&+\frac{k}{mn}f\bigg(\frac{nx+ky}{n+k}\bigg)+\frac{k(n-1)}{mn}f(y)\\&+\bigg(\frac{k(n-1)}{m}+(m-k)k\bigg)\varphi\bigg(\frac{x-y}{n+k}\bigg).
}
Rearranging this inequality, we can conclude that
\Eq{*}{
	f\bigg(\frac{nx+ky}{n+k}\bigg)
	\leq\frac{n(m-k)}{mn-k}f(x)+\frac{k(n-1)}{mn-k}f(y)+kn\frac{(n-1)+m(m-k)}{mn-k}\varphi\bigg(\frac{x-y}{n+k}\bigg).
}
Using that $mn-k=(n-1)(n+k)$, it immediately follows that the inequality \eq{nk} holds if $n+k=m+1$.
\end{proof}

In what follows, we sharpen the inequality \eq{nk} by lowering the error function. Given $\varphi:\tfrac12D_\Delta\to\R$, we define $\varphi^*:\tfrac12D_\Delta\to[-\infty,\infty)$ by
\Eq{*}{
  \varphi^*(u):=\inf_{m\in\N}m^2\varphi\Big(\frac{u}{m}\Big).
}

\Lem{phi}{Let $\varphi:\tfrac12D_\Delta\to\R$. Then $\varphi^*\leq\varphi$ and, for all $k\in\N$, $u\in \tfrac12D_\Delta$,
\Eq{phi*}{
  \varphi^*(u)\leq k^2\varphi^*\Big(\frac{u}{k}\Big).
}
Consequently, $\varphi^{**}:=(\varphi^*)^*=\varphi^*$.
In addition, if $\varphi(0)=0$, then $\varphi^*(0)=0$ and if $\varphi$ is even, then $\varphi^*$ is also even.}

\begin{proof} 
The inequality $\varphi^*\leq\varphi$ is obvious. To verify \eq{phi*}, let $k\in\N$, $u\in \tfrac12D_\Delta$. Then 
\Eq{*}{
  k^2\varphi^*\Big(\frac{u}{k}\Big) 
  =k^2\inf_{m\in\N}m^2\varphi\Big(\frac{u}{km}\Big) 
  =\inf_{m\in\N}(km)^2\varphi\Big(\frac{u}{km}\Big) 
  \geq\inf_{n\in\N}n^2\varphi\Big(\frac{u}{n}\Big) 
  =\varphi^*(u),
}
which proves \eq{phi*}. In view of the inequality \eq{phi*}, for all $u\in \tfrac12D_\Delta$, we have that
\Eq{*}{
  \varphi^*(u)=\inf_{k\in\N}k^2\varphi^*\Big(\frac{u}{k}\Big)=\varphi^{**}(u).
}
Thus, $\varphi^{**}=\varphi^*$. It is obvious that if $\varphi(0)=0$, then $\varphi^*(0)=0$ and if $\varphi$ is even, then $\varphi^*$ is also even.
\end{proof}

\Cor{AQC}{Let $\varphi\in\A(D_\Delta)$. Then, for all $n,k\in\N$ and $x,y\in D$,
\Eq{nk*}{
  f\Big(\dfrac{nx+ky}{n+k}\Big)
  \leq \frac{n}{n+k}f(x)+\frac{k}{n+k}f(y)+
  kn\varphi^*\bigg(\frac{x-y}{n+k}\bigg).
}
Consequently, any $\varphi$-Jensen convex function is also $\varphi^*$-Jensen convex and $\varphi^*\in\A(D_\Delta)$.}

\begin{proof} Let $n,k\in\N$ and $x,y\in D$.
By replacing $(n,k)$ by $(mn,mk)$ in \eq{nk}, for all $m\in\N$, it follows that, 
\Eq{*}{
  f\Big(\dfrac{nx+ky}{n+k}\Big)
  \leq \frac{n}{n+k}f(x)+\frac{k}{n+k}f(y)+
  knm^2\varphi\bigg(\frac{x-y}{m(n+k)}\bigg).
}
Therefore, by the definition of $\varphi^*$, 
\Eq{*}{
  f\Big(\dfrac{nx+ky}{n+k}\Big)
  &\leq \frac{n}{n+k}f(x)+\frac{k}{n+k}f(y)+
  kn\inf_{m\in\N}m^2\varphi\bigg(\frac{x-y}{m(n+k)}\bigg)\\
  &=\frac{n}{n+k}f(x)+\frac{k}{n+k}f(y)+
  kn\varphi^*\bigg(\frac{x-y}{n+k}\bigg),
}
which shows that \eq{nk*} is valid. 

By applying the inequality \eq{nk*} for $n=k=1$, it follows that $f$ is also $\varphi^*$-Jensen convex. 

If $\varphi$ is admissible, then $\varphi$ is even, $\varphi(0)=0$ and there exists $\varphi$-Jensen convex
function, which will also be $\varphi^*$-Jensen convex.
This observation combined with last statements of the previous lemma, yields that $\varphi^*$ is also admissible.
\end{proof}

For $\varphi\in\A(D_\Delta)$, define the \emph{extended error function} $\E[\varphi]:[0,1]\times D_\Delta\to[-\infty,\infty]$ by
\Eq{Psi}{
  \E[\varphi](\lambda,u):=\sup\big\{&f(\lambda x + (1-\lambda)y)
  -\lambda f(x) - (1-\lambda)f(y)\\
  &\mid \mbox{$f:D\to\R$ is $\varphi$-Jensen convex},\ x,y\in D,\,x-y=u\big\}.
}
Since $\varphi\in\A(D_\Delta)$, it follows that $\E[\varphi]$ cannot take the value $-\infty$.
According to this construction, $\E[\varphi]:[0,1]\times D_\Delta\to(-\infty,\infty]$ is the smallest function which, for all $\lambda\in[0,1]$, $x,y\in D$ and for all $\varphi$-Jensen convex functions $f:D\to\R$ satisfies the inequality
\Eq{EE}{
f(\lambda x + (1-\lambda)y)\leq \lambda f(x) + (1-\lambda)f(y) + \E[\varphi](\lambda,x-y).
}

Given $\varphi\in\A(D_\Delta)$, the most important properties of the function $\E[\varphi]$ are contained in the following result. As a consequence of it, we will see that $\E[\varphi]$ has finite values over the set $[0,1]_\Q\times D_\Delta$. 

\Thm{Psi}{Let $\varphi\in\A(D_\Delta)$. Then 
\begin{enumerate}[(i)]
 \item For all $\lambda\in[0,1]$ and $u\in D_\Delta$, $\E[\varphi](\lambda,0)=\E[\varphi](0,u)=\E[\varphi](1,u)=0$.
 \item For all $\lambda\in[0,1]$ and $u\in D_\Delta$, $\E[\varphi](\lambda,u)=\E[\varphi](1-\lambda,-u)$.
 \item For all $n,k\in\N$ and $u\in D_\Delta$,
 \Eq{*}{
   \E[\varphi]\Big(\frac{n}{n+k},u\Big)
   \leq nk\varphi^*\Big(\frac{u}{n+k}\Big).
 }
 \item For all $\lambda\in[0,1]\setminus\Q$ and $u\in D_\Delta\setminus\{0\}$, $\E[\varphi](\lambda,u)=\infty$.
\comment{\item For all $\lambda,\mu\in[0,1]_\Q$ with $1<\lambda+\mu$ and $u\in D_\Delta$,
 \Eq{Psi+}{
  \E[\varphi](\lambda,u)\leq \tfrac{\lambda(1-\lambda)}{\lambda+\mu-1}\E[\varphi]\big(\tfrac{1-\mu}{\lambda},\lambda t\big)+\tfrac{\lambda\mu}{\lambda+\mu-1}\E[\varphi]\big(\tfrac{1-\lambda}{\mu},\mu t\big).
 }}
 \item For all $\lambda,\mu\in[0,1]$ and $u\in D_\Delta$,
\Eq{Psi*}{
   \E[\varphi](\lambda\mu,u)\leq \lambda\E[\varphi]\big(\mu,u\big)+\E[\varphi]\big(\lambda,\mu u\big).
}
\end{enumerate}}

\comment{Let $\mu\in\,]\frac12,1]$ and define $\lambda:=\frac{1-\mu}{\mu}$. Then, from the last inequality, we get
\Eq{*}{
  \E[\varphi](\mu,u)=\E[\varphi](1-\mu,u)
  &\leq \tfrac{1-\mu}{\mu}\E[\varphi]\big(\mu,u\big)
  +\E[\varphi]\big(\tfrac{1-\mu}{\mu},\mu u\big).
}
Interchanging the roles of $\lambda$ and $\mu$ in \eq{Psi*}, we obtain
\Eq{*}{
   \E[\varphi](\mu,u)\leq \tfrac{1-\mu}{\lambda}\E[\varphi]\big(\lambda,t\big)+\E[\varphi]\big(\tfrac{1-\mu}{\lambda},\lambda t\big).
}
Combining this with \eq{Psi*}, we arrive at
\Eq{*}{
  \frac{\lambda+\mu-1}{\lambda\mu}\E[\varphi](\lambda,u)\leq \tfrac{1-\lambda}{\mu}(\E[\varphi]\big(\tfrac{1-\mu}{\lambda},\lambda t\big)+\E[\varphi]\big(\tfrac{1-\lambda}{\mu},\mu t\big)
}}

\begin{proof}
The property (i) directly follows from the definition of $\E[\varphi]$.

To prove (ii), observe that, for all $(\lambda,u)\in[0,1]\times D_\Delta$,
	\Eq{*}{
		\E[\varphi](\lambda,u)&=\sup\big\{f(\lambda x+(1-\lambda)y)-\lambda f(x)-(1-\lambda)f(y)\\
		&\hspace{1.2cm}\mid \mbox{$f:D\to\R$ is $\varphi$-Jensen convex},\ x,y\in D,\,x-y=u\big\}\\
		&=\sup\big\{f((1-\lambda)y+\lambda x)-(1-\lambda)f(y)-\lambda f(x)\\
		&\hspace{1.2cm}\mid \mbox{$f:D\to\R$ is $\varphi$-Jensen convex},\ x,y\in D,\,y-x=-u\big\}\\
		&=\E[\varphi](1-\lambda,-u).
	}
		
In order to show that (iii) is valid, let $n,k\in\N$, $u\in D_\Delta$, and let $c<\E[\varphi]\big(\tfrac{n}{n+k},u\big)$ be arbitrary. Then there exists a $\varphi$-Jensen convex function $f:D\to\R$ and $x,y\in D$ such that $x-y=u$ and
\Eq{*}{
  c<f\Big(\frac{n}{n+k}x +\frac{k}{n+k}y\Big)
  -\frac{n}{n+k}f(x)-\frac{k}{n+k}f(y).
}
According to \cor{AQC}, $f$ is $\varphi^*$-Jensen convex and satisfies \eq{nk*}, therefore,
\Eq{*}{
  c<nk\varphi^*\Big(\frac{x-y}{n+k}\Big)
  =nk\varphi^*\Big(\frac{u}{n+k}\Big).
}
Upon taking the limit $c\to\E[\varphi]\big(\tfrac{n}{n+k},u\big)$, the inequality in condition (iii) follows.

To show that (iv) is valid, let $\lambda\in[0,1]\setminus\Q$ and $u\in D_\Delta\setminus\{0\}$ be fixed.
Then there exists $x,y\in D$ such that $u=x-y$.
Using that $\varphi\in\A(D_\Delta)$, we can choose a $\varphi$-Jensen convex function $f_0:D\to\R$.

Consider $X$ as a vector space over the field $\Q$. Then $u$ and $\lambda u$ are linearly independent elements of this space. Indeed, if for some nonzero pair $(\alpha,\beta)\in\Q^2$, we have that $\alpha u+\beta(\lambda u)=0$, then $\alpha+\beta\lambda=0$, which is impossible by the irrationality of $\lambda$. Then, according to Hamel's theorem, there exists a Hamel base $H$ of $X$ (over the field $\Q$), which contains $u$ and $\lambda u$, furthermore, every function $a:H\to\R$ admits a unique additive extension $A:X\to\R$. Let $\alpha\in\R$ and define now an additive function $A_\alpha:X\to\R$ via the conditions:
\Eq{*}{
  A_\alpha(\lambda u):=\alpha,\qquad A_\alpha(h):=0 \quad(h\in H\setminus\{\lambda u\}).
}
In particular, $A_\alpha(u)=0$.
Then the function $f_\alpha:=f_0+A_\alpha|_D$ is $\varphi$-Jensen convex and
\Eq{*}{
  f_\alpha(\lambda x + (1-\lambda)y)&-\lambda f_\alpha(x) -(1-\lambda)f_\alpha(y)\\
  &=  f_0(\lambda x + (1-\lambda)y)-\lambda f_0(x) -(1-\lambda)f_0(y)\\
  &\quad +  A_\alpha(\lambda x + (1-\lambda)y)-\lambda A_\alpha(x) -(1-\lambda)A_\alpha(y)\\
  &=  f_0(\lambda x + (1-\lambda)y)-\lambda f_0(x) -(1-\lambda)f_0(y)+A_\alpha(\lambda u)-\lambda A_\alpha(u)\\
  &=  f_0(\lambda x + (1-\lambda)y)-\lambda f_0(x) -(1-\lambda)f_0(y)+\alpha.
}
Therefore,
\Eq{*}{
  \E[\varphi]&(\lambda,u)\\
  &\geq\sup\big\{f(\lambda x + (1-\lambda)y)
  -\lambda f(x) - (1-\lambda)f(y)\mid \mbox{$f:D\to\R$ is $\varphi$-Jensen convex}\big\}\\
  &\geq\sup\big\{f_\alpha(\lambda x + (1-\lambda)y)-\lambda f_\alpha(x) -(1-\lambda)f_\alpha(y)\mid \alpha\in\R\big\}\\
  &=\sup\big\{f_0(\lambda x + (1-\lambda)y)-\lambda f_0(x) -(1-\lambda)f_0(y)+\alpha\mid \alpha\in\R\big\}\\
  &=f_0(\lambda x + (1-\lambda)y)-\lambda f_0(x) -(1-\lambda)f_0(y)+\sup\R
  =\infty.
}
	
To prove (v), let $u\in D_\Delta$ and $\lambda,\mu\in[0,1]$ and choose $c<\E[\varphi](\lambda\mu,u)$ arbitrarily. Then there exists a $\varphi$-Jensen convex function $f:D\to\R$ and $x,y\in D$ such that $x-y=u$ and 
\Eq{eps}{
  c<f(\lambda\mu x + (1-\lambda\mu)y)-\lambda\mu f(x) - (1-\lambda\mu)f(y)
}
holds. By the construction of $\E[\varphi]$, we obtain that
\Eq{*}{
  f(\lambda\mu x&+(1-\lambda\mu)y)\\
  &=f\big(\lambda(\mu x+(1-\mu) y)+(1-\lambda)y\big)\\
  &\leq\lambda f(\mu x+(1-\mu) y)+(1-\lambda)f(y)
  +\E[\varphi](\lambda,(\mu x+(1-\mu) y)-y)\\
  &\leq\lambda \Big(\mu f(x)+(1-\mu)f(y)+\E[\varphi]\big(\mu, x-y\big)\Big)+(1-\lambda)f(y)+\E[\varphi](\lambda,\mu u)\\
  &=\lambda\mu f(x) + (1-\lambda\mu)f(y)+\lambda\E[\varphi](\mu, u)+\E[\varphi](\lambda,\mu u).
}
Using inequality \eq{eps}, we obtain that
\Eq{*}{
  c\leq\lambda\E[\varphi](\mu,u)+\E[\varphi](\lambda,\mu u).
}
Upon taking the limit as $c\to\E[\varphi](\lambda\mu,u)$, we get that inequality \eq{Psi*} holds.
\end{proof}

The essential properties of the set $\A(D_\Delta)$ and the functional $\E:\A(D_\Delta)\to(-\infty,\infty]^{[0,1]\times D_\Delta}$ are collected in the following result. 

\Thm{AEF}{Let $D\subseteq X$ be a nonempty convex set. Then $\A(D_\Delta)$ is a convex cone (i.e., $\A(D_\Delta)$ is closed with respect to multiplication by positive scalars and addition), which is also invariant with respect to the mapping $\varphi\mapsto\varphi^*$. Furthermore, for all $t>0$ and $\varphi\in\A(D_\Delta)$,
\Eq{ph}{
  \E[t\varphi]=t\E[\varphi],
}
and, for all $\varphi\in\A(D_\Delta)$,
\Eq{*inv}{
  \E[\varphi]=\E[\varphi^*].
}}

\begin{proof} First we prove that $\A(D_\Delta)$ is closed with respect to multiplication by positive scalars. Let $t>0$, $\varphi\in\A(D_\Delta)$. Then there exists a $\varphi$-Jensen convex function $f:D\to\R$, i.e., \eq{AJC} holds. Multiplying this inequality by $t$ side by side, we can see that $tf$ is a $(t\varphi)$-Jensen convex function. Furthermore, $0\leq(t\varphi)(0)$ and $(t\varphi)$ is also even. Therefore, $t\varphi$  is also an admissible error function.

To show that \eq{ph} holds, let $\lambda\in[0,1]$ and $u\in D_\Delta$, and let $c<\E[\varphi](\lambda,u)$ be arbitrary. Then there exist $x,y\in D$ and a $\varphi$-Jensen convex function $f:D\to\R$ such that $u=x-y$ and 
\Eq{cE}{
  c<f(\lambda x+(1-\lambda)y)-\lambda f(x)-(1-\lambda)f(y).
}
Multiplying this inequality by $t$ side by side, and observing that $tf$ is $(t\varphi)$-Jensen convex, we can see that
\Eq{*}{
  tc<(tf)(\lambda x+(1-\lambda)y)-\lambda (tf)(x)-(1-\lambda)(tf)(y)\leq\E[t\varphi](\lambda,u).
}
Upon taking the limit $c\to\E[\varphi](\lambda,u)$, for all $t>0$, $\varphi\in\A(D_\Delta)$, $\lambda\in[0,1]$, and $u\in D_\Delta$, it follows that
\Eq{ph+}{
  t\E[\varphi](\lambda,u)
  \leq\E[t\varphi](\lambda,u).
}
On the other hand, using this inequality 
\Eq{*}{
  \E[t\varphi](\lambda,u)
  =t\cdot\frac{1}{t}\E[t\varphi](\lambda,u)
  \leq t\E[\tfrac{1}{t}\cdot (t\varphi)](\lambda,u)
  =t\E[\varphi](\lambda,u),
}
which shows that \eq{ph+} holds with equality.
This shows that \eq{ph} is also valid.

Next we point out that $\A(D_\Delta)$ is closed with respect to addition. Let $\varphi,\psi \in\A(D_\Delta)$. Then there exist a $\varphi$-Jensen convex function $f:D\to\R$ and a $\psi$-Jensen convex function $g:D\to\R$. It is immediate to see that $f+g$ is a $(\varphi+\psi)$-Jensen convex function. Furthermore, $0\leq(\varphi+\psi)(0)$ and $\varphi+\psi$ is also even. Therefore, $\varphi+\psi$  is also an admissible error function.

In view of the last statement of \cor{AQC}, if $\varphi\in\A(D_\Delta)$, then $\varphi^*\in\A(D_\Delta)$. The inequality $\varphi^*\leq\varphi$ implies that $\E[\varphi^*]\leq\E[\varphi]$. To prove the reversed inequality, let $\lambda\in[0,1]$ and $u\in D_\Delta$, and let $c<\E[\varphi](\lambda,u)$ be arbitrary. Then there exist $x,y\in D$ and a $\varphi$-Jensen convex function $f:D\to\R$ such that $u=x-y$ and the inequality \eq{cE} holds. By \cor{AQC}, we have that $f$ is also $\varphi^*$-Jensen convex, hence
\Eq{*}{
  f(\lambda x+(1-\lambda)y)-\lambda f(x)-(1-\lambda)f(y)\leq\E[\varphi^*](\lambda,u).
}
Therefore, $c<\E[\varphi^*](\lambda,u)$. Upon taking the limit $c\to\E[\varphi](\lambda,u)$, it follows that the inequality $\E[\varphi](\lambda,u)\leq\E[\varphi^*](\lambda,u)$ is valid for all $\lambda\in[0,1]$ and $u\in D_\Delta$. Consequently, $\E[\varphi]\leq\E[\varphi^*]$. Thus, the equality
\eq{*inv} is proved.
\end{proof}

In the following result we establish an upper estimate for $\E[\varphi](\lambda,u)$ obtained in assertion (iii) of \thm{Psi} in terms of a Takagi-type function related to the map $t\mapsto \varphi^*(tu)$, where $u\in D_\Delta$ is fixed.

\Thm{MP+}{Let $\varphi\in\A(D_\Delta)$. Then, for all $\lambda\in[0,1]_\Q$ and $u\in D_\Delta$,
\Eq{MP+}{
  \E[\varphi](\lambda,u)
  \leq\sum_{k=0}^\infty
  \frac{1}{2^k}\varphi^*\big(d_\Z(2^k\lambda)u\big),
}
furthermore, the series on the right hand side of \eq{MP+} is absolutely convergent and if $\lambda=\frac{n}{n+m}$ for some $n,m\in\N$, then, for all $u\in D_\Delta$,
\Eq{MP++}{
  \sum_{k=0}^\infty
  \frac{1}{2^k}\varphi^*\big(d_\Z(2^k\lambda)u\big)
  \leq nm\varphi^*\Big(\frac{u}{n+m}\Big).
}}

\begin{proof} 
First, we verify that the series on the right hand side of \eq{MP+} is absolutely convergent for all $\lambda\in[0,1]_\Q$ and $u\in D_\Delta$. For $\lambda\in[0,1]$, define the set $S(\lambda)\subseteq \big[1,\frac12\big]$ by
\Eq{*}{
  S(\lambda):=\big\{d_\Z(2^k\lambda)\colon k\in\N\cup\{0\}\big\}.
}
We claim that this set is finite for all $\lambda\in[0,1]_\Q$. If $\lambda\in\{0,1\}$, then $d_\Z(2^k\lambda)=0$ for all $k\in\N\cup\{0\}$, whence we get that $S(0)=S(1)=\{0\}$. Thus, we may assume that $\lambda\in(0,1)_\Q$. In this case $\lambda$ can be written uniquely in the form $\lambda=2^\ell\frac{m}{n}$, where $\ell\in\Z$ and $n,m$ are odd coprime natural numbers. If $k\geq-\ell$ the denominator of $2^k\lambda=2^{k+\ell}\frac{m}{n}$ is equal to $n$. Therefore, the distance of $2^k\lambda$ from $\Z$ is a rational number whose denominator equals $n$, i.e.,
\Eq{*}{
  \big\{d_\Z(2^k\lambda)\colon k\in\N\cup\{0\},\,k\geq-\ell\}\subseteq\big\{\tfrac{i}{n}\colon i\in\big\{0,1,\dots,\big\lfloor\tfrac{n}{2}\big\rfloor\big\}\big\}.
}
Therefore,
\Eq{*}{
  S(\lambda)
  &=\big\{d_\Z(2^k\lambda)\colon k\in\N\cup\{0\},\,k<-\ell\}\cup\big\{d_\Z(2^k\lambda)\colon k\in\N\cup\{0\},\,k\geq-\ell\}\\
  &\subseteq\big\{d_\Z(2^k\lambda)\colon k\in\N\cup\{0\},\,k<-\ell\}\cup\big\{\tfrac{i}{n}\colon i\in\big\{0,1,\dots,\big\lfloor\tfrac{n}{2}\big\rfloor\big\}\big\},
}
which shows that $S(\lambda)$ is covered by union of two finite sets. Thus $S(\lambda)$ is finite as well. 

Now we are able to show that, $\lambda\in[0,1]_\Q$ and $u\in D_\Delta$, the series on the right hand side of \eq{MP+} is absolutely convergent. Indeed,
\Eq{*}{
  \sum_{k=0}^\infty
  \frac{1}{2^k}\big|\varphi^*\big(d_\Z(2^k\lambda)u\big)\big|
  \leq \sum_{k=0}^\infty
  \frac{1}{2^k}\max_{s\in S(\lambda)}\big|\varphi^*(s u)\big|
  =2\max_{s\in S(\lambda)}\big|\varphi^*(s u)\big|<\infty.
}

In what follows, we prove by induction with respect to $n\in\N$ that, for all $\lambda\in[0,1]_\Q$ and $u\in D_\Delta$,
\Eq{n}{
  \E[\varphi](\lambda,u)
  \leq \sum_{k=0}^{n-1}
  \frac{1}{2^k}\varphi^*\big(d_\Z(2^k\lambda)u\big)
  +\frac{\pmb\langle\lambda\pmb\rangle^2}{2^{n+2}} \max_{1\leq m\leq\pmb\langle\lambda\pmb\rangle}\left|\varphi^*\bigg(\frac{u}{m}\bigg)\right|.
}
If $\lambda=\frac{n}{n+k}$, where $n,k\in\N$ are coprime numbers, then $\pmb\langle\lambda\pmb\rangle=n+k$.
By the inequality $nk\leq\frac{1}{4}(n+k)^2$ and assertion (iii) of \thm{Psi}, it follows that
\Eq{*}{
  \E[\varphi](\lambda,u)
  &\leq nk \varphi^*\bigg(\frac{u}{n+k}\bigg)
  \leq nk \left|\varphi^*\bigg(\frac{u}{n+k}\bigg)\right|\\
  &\leq \frac{(n+k)^2}{4} \left|\varphi^*\bigg(\frac{u}{n+k}\bigg)\right|
  =\frac{\pmb\langle\lambda\pmb\rangle^2}{2^2} \left|\varphi^*\bigg(\frac{u}{\pmb\langle\lambda\pmb\rangle}\bigg)\right|
  \leq\frac{\pmb\langle\lambda\pmb\rangle^2}{2^2} \max_{1\leq m\leq\pmb\langle\lambda\pmb\rangle}\left|\varphi^*\bigg(\frac{u}{m}\bigg)\right|,
}
which proves the inequality \eq{n} for $n=0$. 

Assume that we have proved \eq{n} for some $n\geq0$. First let $\lambda\in[0,\frac12]_\Q$ and let $u\in D_\Delta$ be fixed. By assertion (v) of \thm{Psi}, we have
\Eq{1.0}{
  \E[\varphi](\lambda,u)
  =\E[\varphi](\tfrac12\cdot(2\lambda),u)
  \leq \tfrac12\E[\varphi]\big(2\lambda,u\big)+\E[\varphi]\big(\tfrac12,2\lambda u\big).
}
By assertion (iii) of \thm{Psi}, we have that
\Eq{1.1}{
  \E[\varphi]\big(\tfrac12,2\lambda u\big)\leq\varphi^*\big(\lambda u\big)
  =\varphi^*\big(d_\Z(\lambda) u\big).
}
In view of the inductive hypothesis, we can obtain
\Eq{1.2}{
  \frac12\E[\varphi]\big(2\lambda,u\big)
 &\leq \frac12\bigg(\sum_{k=0}^{n-1}
  \frac{1}{2^k}\varphi^*\big(d_\Z(2^k(2\lambda))u\big)
  +\frac{\pmb\langle2\lambda\pmb\rangle^2}{2^{n+2}}\max_{1\leq m\leq\pmb\langle2\lambda\pmb\rangle}\left|\varphi^*\bigg(\frac{u}{m}\bigg)\right|\bigg)\\
  &=\sum_{k=1}^{n}
  \frac{1}{2^{k}}\varphi^*\big(d_\Z(2^{k}\lambda))u\big)
  +\frac{\pmb\langle2\lambda\pmb\rangle^2}{2^{n+3}}\max_{1\leq m\leq\pmb\langle2\lambda\pmb\rangle}\left|\varphi^*\bigg(\frac{u}{m}\bigg)\right|\\
  &\leq\sum_{k=1}^{n}
  \frac{1}{2^{k}}\varphi^*\big(d_\Z(2^{k}\lambda))u\big)
  +\frac{\pmb\langle\lambda\pmb\rangle^2}{2^{n+3}}\max_{1\leq m\leq\pmb\langle\lambda\pmb\rangle}\left|\varphi^*\bigg(\frac{u}{m}\bigg)\right|.
}
Applying the inequalities \eq{1.1} and \eq{1.2}, the inequality \eq{1.0} implies that \eq{n} holds for all $\lambda\in[0,\frac12]_\Q$, $u\in D_\Delta$ and for $n+1$ (instead of $n$).

Now assume that $\lambda\in[\frac12,1]_\Q$ and let $u\in D_\Delta$ be fixed. Then $1-\lambda\in[0,\frac12]_\Q$. Using assertion (ii) of \thm{Psi}, and then what we have proved in the previous case, we get
\Eq{*}{
  \E[\varphi](\lambda,u)
  &=  \E[\varphi](1-\lambda,-u)\\
  &\leq\sum_{k=0}^{n-1}
  \frac{1}{2^k}\varphi^*\big(d_\Z(2^k(1-\lambda))(-u)\big)
  +\frac{\pmb\langle1-\lambda\pmb\rangle^2}{2^{n+2}} \max_{1\leq m\leq\pmb\langle1-\lambda\pmb\rangle}\left|\varphi^*\bigg(\frac{-u}{m}\bigg)\right|\\
  &\leq\sum_{k=0}^{n-1}
  \frac{1}{2^k}\varphi^*\big(d_\Z(2^k\lambda)u\big)
  +\frac{\pmb\langle\lambda\pmb\rangle^2}{2^{n+2}} \max_{1\leq m\leq\pmb\langle\lambda\pmb\rangle}\left|\varphi^*\bigg(\frac{u}{m}\bigg)\right|.
}
Thus, \eq{n} also holds for all $\lambda\in[\frac12,1]$, $u\in D_\Delta$ and for $n+1$ (instead of $n$).

This completes the proof of inequality \eq{n} on the indicated domain. Upon taking the limit $n\to\infty$, we can conclude that the asserted inequality is valid for all $\lambda\in[0,1]_\Q$, $u\in D_\Delta$.

Finally, we show that the inequality \eq{MP++} is also valid. Assume that $\lambda=\frac{n}{n+m}$ for some $n,m\in\N$. Then $(n+m)d_\Z(2^k\lambda)$ is an integer number for all $k\in\N$. Therefore, due to the inequality \eq{phi*}, we have that
\Eq{*}{
  \varphi^*\big(d_\Z(2^k\lambda)u\big)
  =\varphi^*\Big((n+m)d_\Z(2^k\lambda)\frac{u}{n+m}\Big)
  \leq (n+m)^2 d_\Z(2^k\lambda)^2\varphi^*\Big(\frac{u}{n+m}\Big).
}
Therefore, using \cor{psi}, we get
\Eq{*}{
  \sum_{k=0}^\infty \frac{1}{2^k}\varphi^*\big(d_\Z(2^k\lambda)u\big)
  &\leq \sum_{k=0}^\infty \frac{(n+m)^2 d_\Z(2^k\lambda)^2}{2^k}\varphi^*\Big(\frac{u}{n+m}\Big)\\
  &\leq \bigg(\sum_{k=0}^\infty \frac{ d_\Z(2^k\lambda)^2}{2^k}\bigg)(n+m)^2\varphi^*\Big(\frac{u}{n+m}\Big)\\
  &=\lambda(1-\lambda)(n+m)^2\varphi^*\Big(\frac{u}{n+m}\Big)
  =nm\varphi^*\Big(\frac{u}{n+m}\Big).
}
This completes the proof of the inequality \eq{MP++}.
\end{proof}

In our last result, we establish a closed formula for the computation of the Takagi-type function given by the infinite series on the right hand side of the inequality \eq{MP+}.

\Thm{Last}{Let $\varphi\in\A(D_\Delta)$ and $\lambda:=\frac{m}{2^jn}$, where $m,j\geq0$ are integers, $n\geq1$ is an odd integer, $m\leq 2^{j-1}n$, furthermore, $m$ and $2^jn$ are coprime numbers. Denote $\ell:=\left\lfloor\frac12\phi(n)\right\rfloor$. Then, for all $u\in D_\Delta$,
\Eq{Last}{
  \sum_{k=0}^\infty \frac{1}{2^k}\varphi^*\big(d_\Z(2^k\lambda)u\big)
  =\sum_{k=0}^{j-1} \frac{1}{2^k}\varphi^*\big(d_\Z(2^k\lambda)u\big)
  +\sum_{k=j}^{j+\ell-1}\frac{1}{2^k-2^{k-\ell}}\varphi^*\big(d_\Z(2^k\lambda)u\big).
}}

\begin{proof}
If $n=1$, then, for $k\geq j$, the product $2^k\lambda$ is an integer. Thus, $d_\Z(2^k\lambda)=0$, which implies that $\varphi^*\big(d_\Z(2^k\lambda)u\big)=\varphi^*(0)=0$. Therefore,
\Eq{*}{
  \sum_{k=0}^\infty \frac{1}{2^k}\varphi^*\big(d_\Z(2^k\lambda)u\big)
  =\sum_{k=0}^{j-1} \frac{1}{2^k}\varphi^*\big(d_\Z(2^k\lambda)u\big).
}
This shows that \eq{Last} holds in the case $n=1$.

In the rest of the proof, we assume that $n\geq3$. Then $\ell=\frac12\phi(n)$ (since $\phi(n)$ is even for $n\geq2$). Based on \prp{kl}, for $i,k\geq0$, we have that 
\Eq{*}{
  d_\Z(\tfrac{2^{k}\cdot m}{n})
  =d_\Z(2^{i\ell}\cdot\tfrac{2^{k}\cdot m}{n}).
}
From this equality it follows that, for $i\geq0$ and $k\geq j$, 
\Eq{*}{
  d_\Z(2^{k}\lambda)
  =d_\Z(2^{k+i\ell}\lambda).
}
In view of this equality, we obtain
\Eq{*}{
  \sum_{k=0}^\infty \frac{1}{2^k}\varphi^*\big(d_\Z(2^k\lambda)u\big)
  &=\sum_{k=0}^{j-1} \frac{1}{2^k}\varphi^*\big(d_\Z(2^k\lambda)u\big)
  +\sum_{k=j}^{\infty} \frac{1}{2^k}\varphi^*\big(d_\Z(2^k\lambda)u\big)\\
  &=\sum_{k=0}^{j-1} \frac{1}{2^k}\varphi^*\big(d_\Z(2^k\lambda)u\big)
  +\sum_{k=j}^{j+\ell-1}\sum_{i=0}^\infty \frac{1}{2^{k+i\ell}}\varphi^*\big(d_\Z(2^{k+i\ell}\lambda)u\big)\\
  &=\sum_{k=0}^{j-1} \frac{1}{2^k}\varphi^*\big(d_\Z(2^k\lambda)u\big)
  +\sum_{k=j}^{j+\ell-1}\bigg(\sum_{i=0}^\infty \frac{1}{2^{k+i\ell}}\bigg)\varphi^*\big(d_\Z(2^{k}\lambda)u\big)\\
  &=\sum_{k=0}^{j-1} \frac{1}{2^k}\varphi^*\big(d_\Z(2^k\lambda)u\big)
  +\sum_{k=j}^{j+\ell-1}\frac{1}{2^{k}-2^{k-\ell}}\varphi^*\big(d_\Z(2^{k}\lambda)u\big),
}
which completes the proof of the formula \eq{Last}.
\end{proof}

In what follows, we  present a few particular cases of the combination of inequality \eq{MP+} and the formula \eq{Last}. Let $u\in D_\Delta$ be fixed.

If $\lambda=\frac12$, then $j=n=m=1$ and hence $\ell=0$. Thus, the second sum in \eq{Last} is empty and we get,
\Eq{*}{
  \E[\varphi](\tfrac{1}{2},u)
  \leq\sum_{k=0}^0\frac{1}{2^k}\varphi^*\big(2^{k-1}u\big)
  =\varphi^*\big(\tfrac12 u\big).
}

If $\lambda=\frac13$ or $\lambda=\frac23$, then $j=0$, $n=3$, $m\in\{1,2\}$ and $\ell=1$. Therefore, the first sum in \eq{Last} is empty and we get,
\Eq{*}{
  \E[\varphi](\tfrac{i}{3},u)
  \leq \sum_{k=0}^{0}\frac{1}{2^k-2^{k-1}}\varphi^*\big(d_\Z(2^k\lambda)u\big)
  = 2\varphi^*\big(\tfrac{1}{3}u\big)
  \qquad (i\in\{1,2\}).
}
It is  visible that, in these two cases, the inequality \eq{MP+} is equivalent to the assertion (iii) of \thm{Psi}. On the other hand, in the particular cases below, we obtain sharper upper estimates.

If $\lambda=\frac14$ or $\lambda=\frac34$, then $j=2$, $n=1$, $m\in\{1,3\}$, and $\ell=0$. In this case, the second sum in \eq{Last} is empty and we get,
\Eq{*}{
  \E[\varphi](\tfrac{i}{4},u)
  \leq\sum_{k=0}^1\frac{1}{2^k}\varphi^*\big(d_\Z(2^k\lambda)u\big)
  =\varphi^*\big(\tfrac{1}{4}u\big)+\frac12\varphi^*\big(\tfrac{1}{2}u\big)
  \qquad (i\in\{1,3\}).
}

If $\lambda=\frac{1}{5}$ or $\lambda=\frac{4}{5}$, then $j=0$, $n=5$, $m\in\{1,2,3,4\}$, and $\ell=2$. In this case, the first sum in \eq{Last} is empty and we obtain,
\Eq{*}{
	\E[\varphi](\tfrac{i}{5},u)
	&\leq
	\sum_{k=0}^{1}\frac{1}{2^k-2^{k-2}}\varphi^*\big(d_\Z(2^k\lambda)u\big)
	=\begin{cases}
    \dfrac{4}{3}\varphi^*\big(\tfrac{1}{5}u\big)
	+\dfrac{2}{3}\varphi^*\big(\tfrac{2}{5}u\big)
	&\mbox{if } i\in\{1,4\},\\[4mm]
    \dfrac{4}{3}\varphi^*\big(\tfrac{2}{5}u\big)
	+\dfrac{2}{3}\varphi^*\big(\tfrac{1}{5}u\big)
	&\mbox{if } i\in\{2,3\}.
	\end{cases}
}

If $\lambda=\frac{1}{6}$ or $\lambda=\frac{5}{6}$, then $j=1$, $n=3$, $m\in\{1,5\}$, and $\ell=1$. In this case, the both sums in \eq{Last} are nonempty and, for $i\in\{1,5\}$, we obtain,
\Eq{*}{
  \E[\varphi](\tfrac{i}{6},u)
  \leq \sum_{k=0}^{0} \frac{1}{2^k}\varphi^*\big(d_\Z(2^k\lambda)u\big)
  +\sum_{k=1}^{1}\frac{1}{2^k-2^{k-1}}\varphi^*\big(d_\Z(2^k\lambda)u\big)
  =\varphi^*\big(\tfrac{1}{6}u\big)
    +\varphi^*\big(\tfrac{1}{3}u\big).
}


\def\MR#1{}
\providecommand{\MRhref}[2]{%
  \href{http://www.ams.org/mathscinet-getitem?mr=#1}{#2}
}
\providecommand{\href}[2]{#2}

\end{document}